\documentclass[12pt,oneside,reqno]{amsart}
\usepackage{graphicx}
\usepackage{mathrsfs}
\usepackage{stmaryrd}
\usepackage{amsfonts}
\usepackage{enumerate,amsmath,amssymb,amsthm}
\newcommand{\be}{\begin{eqnarray}}
\newcommand{\ee}{\end{eqnarray}}
\newcommand{\ce}{\begin{eqnarray*}}
\newcommand{\de}{\end{eqnarray*}}
\newtheorem{theorem}{Theorem}[section]
\newtheorem{lemma}[theorem]{Lemma}
\newtheorem{remark}[theorem]{Remark}
\newtheorem{definition}[theorem]{Definition}
\newtheorem{proposition}[theorem]{Proposition}

\newtheorem{corollary}[theorem]{Corollary}
\def\eps{\varepsilon}

\def\Om{\Omega}

\def\p{\partial}

\def\[{{\Big[}}
\def\]{{\Big]}}
\def\<{{\langle}}
\def\>{{\rangle}}
\def\({{\Big(}}
\def\){{\Big)}}

\def\dif{{\mathord{{\rm d}}}}

\def\no{\nonumber}
\def\bt{\begin{theorem}}
\def\et{\end{theorem}}
\def\bl{\begin{lemma}}
\def\el{\end{lemma}}
\def\br{\begin{remark}}
\def\er{\end{remark}}
\def\bd{\begin{definition}}
\def\ed{\end{definition}}
\def\bp{\begin{proposition}}
\def\ep{\end{proposition}}
\def\bc{\begin{corollary}}
\def\ec{\end{corollary}}
\def\cA{{\mathcal A}}
\def\cB{{\mathcal B}}

\def\cF{{\mathcal F}}

\def\cM{{\mathcal M}}

\def\cO{{\mathcal O}}

\def\mB{{\mathbb B}}
\def\mC{{\mathbb C}}

\def\mE{{\mathbb E}}

\def\mH{{\mathbb H}}

\def\mK{{\mathbb K}}
\def\mL{{\mathbb L}}

\def\mN{{\mathbb N}}

\def\mQ{{\mathbb Q}}
\def\mR{{\mathbb R}}
\def\mS{{\mathbb S}}

\def\mU{{\mathbb U}}

\def\mX{{\mathbb X}}
\def\mY{{\mathbb Y}}

\def\sB{{\mathscr B}}

\def\geq{\geqslant}
\def\leq{\leqslant}

\allowdisplaybreaks

\begin{document}

\title{Freidlin-Wentzell's Large Deviations for Stochastic Evolution Equations}

\author{Jiagang Ren$^{1}$, Xicheng Zhang$^{2,3}$}

\date{}
\dedicatory{
$^1$School of Mathematics and Computational Science\\
Zhongshan University, Guangzhou,
Guangdong 510275, P.R.China\\
$^2$School of Mathematics and Statistics\\
The University of New South Wales, Sydney, 2052, Australia,\\
$^3$Department of Mathematics,
Huazhong University of Science and Technology,\\
Wuhan, Hubei 430074, P.R.China\\
emails: J. Ren: mcdm06@zsu.edu.cn;
X. Zhang: XichengZhang@gmail.com.}

\begin{abstract}
We prove a Freidlin-Wentzell large deviation principle for general
stochastic evolution equations with small perturbation multiplicative noises.
In particular, our general result can be used to deal with a large class of
quasi linear stochastic partial differential equations,
such as stochastic porous medium equations and
stochastic reaction diffusion equations with polynomial growth zero order term
and $p$-Laplacian second order term.

\end{abstract}


\maketitle
\allowdisplaybreaks
\section{Introduction}

Since the work of Freidlin and Wentzell \cite{Fr-We}, the theory of
small perturbation large deviations for stochastic differential equations(SDE)
has been extensively developed(cf. \cite{Az,St}, etc.).
In classical method, to establish such a large deviation principle(LDP)
for SDEs, one needs to discretize the time variable and then
prove various necessary exponential continuity and tightness
for stochastic dynamical systems in different spaces by using comparison principle.
However, such verifications would become rather complicated and even impossible
in some cases for infinite stochastic partial differential equations
with multiplicative noises.

Recently, Dupuis and Ellis \cite{de} systematically developed
a weak convergence approach to the theory of large deviation.
The core idea is to prove some variational representation formula
about the Laplace transform of bounded continuous functionals,
which will lead to proving an equivalent Laplace principle with LDP.
In particular, for Brownian functionals, an elegant variational representation
formula has been established by Bou\'e-Dupuis \cite{bd}
and Budhiraja-Dupuis \cite{bd0}. A simplified proof is given by the second named author
\cite{Zh2}. This variational representation has been proved to be very effective
for various finite dimensional stochastic dynamical system
with irregular coefficients(cf. \cite{bde,Rz1,Rz2}, etc.).
One of the main advantages of this argument is that one only needs to
make some necessary moment estimates.  This can be seen completely from the
present paper that it also works very well for infinite dimensional
stochastic dynamical systems.

In the past two decades, there are  numerous results about
the LDP for stochastic partial differential equations(SPDE)
(cf. \cite{So,DaZa,Pe,Ka-Xi,Ce-Ro,Fe-Ku,Bu-Du-Ma}, etc.).
All these results are concentrated on semi-linear SPDEs, i.e.,
the second order term is linear, and their proofs, except \cite{Fe-Ku,Bu-Du-Ma},
are mainly based on the classical exponential tightness method.
In \cite{Fe-Ku}, the approach for LDP is based on nonlinear semigroup and infinite dimensional
Hamilton-Jacobi equations. The approach in \cite{Bu-Du-Ma} is based on the variational
representation.
Recently, R\"ockner-Wang-Wu \cite{Ro-Wa-Wu} proved an LDP for stochastic porous medium equation
with additive noise by using the classical comparison principle.
It should be pointed out that the equation of this type has a non-linear and degenerated
second order term. Since additive noise was considered in \cite{Ro-Wa-Wu},
they can discretize time and prove some necessary estimates.
It seems difficult to extend their result to the multiplicative noise case
by using the classical method.

On the other hand, the existence and uniqueness of SPDEs have already been
studied in various literatures
prior to  LDP for SPDEs(cf. \cite{Pa1, Kr-Ro,Ro,DaZa,Go-Mi,Zh1}, etc.).
In the theory of SPDEs, there exist two main tools: semigroup method and
variation method(or monotone method).
One of the merits of semigroup method is that the noise can take values
in a larger space(cf. \cite{DaZa}).
But, it can only deal with semi-linear SPDEs. The variation method combined with
Galerkin's approximation is usually used in the framework of evolution triple(cf. \cite{Kr-Ro,Zh1}).
Thus, as in the deterministic case(cf. \cite{Sh}), it can tackle a large class of SPDEs.
But, the diffusion coefficients need to be in the space of Hilbert-Schdmit operators.

Our aim in the present paper is to prove a Freidlin-Wentzell's large deviation
for stochastic evolution equations in the evolution triple case by using the weak convergence
approach as done in \cite{Bu-Du-Ma}. Thus, the main point
is to prove the tightness of some control stochastic evolution equations.
This will be realized by making some moment estimates in suitable space(see Lemma \ref{Le2} below)
and then using the general tightness criterion for stochastic processes(see Lemma \ref{Le4} below).
Moreover, in order to treat the SPDEs with polynomial growth,
we will work in the framework of \cite{Zh1}, which is a little different from \cite{Kr-Ro}.
Compared with the well-known results, our proof is succinct, and
we believe that our method can be adapted to
some other non-linear stochastic equations such as
stochastic Navier-Stokes equation.

This paper is organized as follows: In Section 2, we shall give our framework and
recall an abstract criterion for Laplace principle due to Budhiraja-Dupuis \cite{bd0},
as well as an existence and uniqueness result for stochastic evolution equation essentially
due to Krylov-Rozovskii \cite{Kr-Ro}. In Section 3,
we first prove a Laplace principle for stochastic evolution equation(see Theorem \ref{Th1} below)
without any compact embedding requirement.
In order to prove the corresponding rate function is good, we need an extra compact assumption
(see Lemma \ref{Th2} below).
Lastly, in Section 4 we give three applications.

\section{Framework and Preliminaries}

Let $\mX$ be a reflexive and separable Banach space, which is densely
and continuously injected in a separable Hilbert space
$\mH$. Identifying $\mH$ with its dual we get
$$
\mX\subset\mH\simeq\mH^*\subset\mX^*,
$$
where the star `$^*$' denotes the dual spaces.

Assume that the norm in $\mX$ is given by
$$
\|x\|_\mX:=\|x\|_{1,\mX}+\|x\|_{2,\mX},\quad x\in\mX.
$$
Denote by $\mX_i$, $i=1,2$ the completions of $\mX$ with respect to the norms $\|\cdot\|_{i,\mX}
=:\|\cdot\|_{\mX_i}$. Then $\mX=\mX_1\cap\mX_2$.
Let us also assume that both spaces are reflexive and embedded in $\mH$.
Thus, we get two triples:
$$
\mX_1\subset\mH\simeq\mH^*\subset\mX^*_1, \ \ \mX_2\subset\mH\simeq\mH^*\subset\mX^*_2.
$$
Noticing that $\mX^*_1$ and $\mX^*_2$ can be thought as subspaces of $\mX^*$, one may define
a  Banach space $\mY:=\mX^*_1+\mX^*_2\subset\mX^*$ as follows: $f\in\mY$ if and only if
$f=f_1+f_2$, $f_i\in\mX^*_i, i=1,2$ and the norm of $f$ is defined by
$$
\|f\|_{\mY}=\inf_{f=f_1+f_2}(\|f_1\|_{\mX^*_1}+\|f_2\|_{\mX^*_2}).
$$
In the following, the dual pairs of $(\mX,\mX^*)$ and $(\mX_i,\mX^*_i), i=1,2$ are denoted respectively by
$$
[\cdot,\cdot]_{\mX},\quad[\cdot,\cdot]_{\mX_i}, \ \ \ i=1,2.
$$
Then, for any $x\in\mX$ and $f=f_1+f_2\in\mY\subset\mX^*$,
$$
[x,f]_{\mX}=[x,f_1]_{\mX_1}+[x,f_2]_{\mX_2}.
$$
We remark that if $f\in\mH$ and $x\in\mX$, then
\be
[x,f]_{\mX}=[x,f]_{\mX_1}=[x,f]_{\mX_2}=\<x,f\>_\mH,\label{PP4}
\ee
where $\<\cdot,\cdot\>_\mH$ stands for the inner product in $\mH$.

Let $(\Om,\cF,(\cF_t)_{t\geq 0}, P)$ be a complete separable filtration probability space, and
$Q$ a nonnegative definite and symmetric trace operator defined on another separable
Hilbert space $\mU$. A $Q$-Wiener process $\{W(t),t\geq0\}$ defined on
$(\Om,\cF,P)$  is given and assumed to be adapted to $(\cF_t)_{t\geq 0}$(cf. \cite{DaZa}).
Set $\mU_Q:=Q^{1/2}(\mU)$ and let $L_2(\mU_Q,\mH)$ denote the Hilbert space consisting of
all Hilbert-Schmidt operators from $\mU_Q$ to $\mH$, where
the inner product is denoted by $\<\cdot,\cdot,\>_{L_2(\mU_Q,\mH)}$, and the norm
by $\|\cdot\|_{L_2(\mU_Q,\mH)}$.

In the following, we will work in the finite time interval $[0,T]$.
For a Banach space $\mB$ we shall denote by
$\mC_T(\mB)$ the continuous functions space from $[0,T]$
to $\mB$, which is endowed with the uniform norm.
Define
$$
\mL_Q:=\left\{h=\int^\cdot_0\dot h(s)\dif s: ~~\dot h\in L^2(0,T;\mU_Q)\right\}
$$
with the norm
$$
\|h\|_{\mL_Q}:=\left(\int^1_0\|\dot h(s)\|_{\mL_Q}^2\dif s\right)^{1/2},
$$
where the dot denotes the generalized derivative.
Let $\mu_Q$ be the law of the $Q$-Wiener process $W$ in $\mC_T(\mU)$. Then
$$
(\mC_T(\mU), \mL_Q,\mu_Q)
$$
forms an abstract Wiener space.

For $N>0$ we set $D_N:=\{h\in \mL_Q: \|h\|_{\mL_Q}\leq N\}$.
Then $D_N$ is metrizable as a compact Polish space with respect to the weak
topology in $\mL_Q$. Let $\cA_N$ denote all continuous and $\cF_t$-adapted process $h$ from $[0,T]$ to $\mU_Q$
such that for almost all $\omega$, $h(\cdot,\omega)\in D_N$, i.e.,
\be
\int^T_0\|\dot h(s,\omega)\|_{\mU_Q}^2\dif s\leq N.\label{Op2}
\ee

Let $\mS$ be a Polish space. A function $I: \mS\to[0,\infty]$ is given.
\bd
The function $I$  is called a rate function if $I$ is lower semicontinuous.
The function $I$ is called a good rate function if for every $a<\infty$,
$\{f\in\mS: I(f)\leq a\}$ is compact.
\ed

Let $Z^\eps: \mC_T(\mU)\to\mS$
be a family of measurable mappings. We assume that\\
\begin{enumerate}[{\bf (Hypothesis)}:]
\item
There is a measurable map $Z^0: \mL_Q\mapsto \mS$ such that
for any $N>0$, if a family $\{h^\eps\}\subset\cA_N$(as random variables in $D_N$)
converges in distribution
to a $v\in \cA_N$, then  for some subsequence $\eps_k$,
$Z^{\eps_k}(\cdot+\frac{h^{\eps_k}(\cdot)}{\sqrt{\eps_k}})$ converges
in distribution to $Z^0(v)$ in $\mS$.
\end{enumerate}
For each $f\in\mS$, define
\be
I(f):=\frac{1}{2}\inf_{\{h\in\mL_Q:~f=Z^0(h)\}}\|h\|^2_{\mL_Q},\label{ra}
\ee
where $\inf\emptyset=\infty$ by convention.

We recall the following result due to \cite{bd, bd0}(see also \cite[Theorem 4.4]{Zh2}).
\bt\label{Th2}
$\{Z^\eps,\eps\in(0,1)\}$ satisfies the Laplace principle with
the rate function $I(f)$ given by (\ref{ra}). That is, for each real bounded continuous
function  $g$ on $\mS$:
\be
\lim_{\eps\rightarrow 0}\eps\log\mE\left(\exp\left[-\frac{g(Z^\eps)}{\eps}\right]\right)
=-\inf_{f\in\mS}\{g(f)+I(f)\}.\label{La}
\ee
\et

\br
If $I$ in (\ref{La}) is not lower semicontinuous, then the regularization of
$I$
$$
\tilde I(f):=\lim_{\eps\downarrow 0}\inf_{f'\in B_\eps(f)}I(f)
$$
still satisfies (\ref{La}), where $B_\eps(f)$ is the ball in $\mS$ with center
$f$ and radius $\eps$. Moreover, if $I$ is a good rate function,
then the Laplace principle is equivalent to
the large deviation principle(cf. \cite[Theorem 1.2.3]{de}).
\er

We now introduce three evolution  operators used in the present paper(cf. \cite{Zh1}):
$$
A_i: [0,T]\times\mX_i\to\mX^*_i\in \cB([0,T])\times\cB(\mX_i)/\cB(\mX^*_i),\ \ i=1,2,
$$
and
$$
B: [0,T]\times\mH\to L_2(\mU_Q,\mH)\in \cB([0,T])\times\cB(\mH)/\cB(L_2(\mU_Q,\mH)).
$$
In the following, for the sake of simplicity, we write
$$
A=A_1+A_2\in\mY\subset\mX^*,
$$
and assume  throughout this paper that
\begin{enumerate}[{\bf (H1)}]
\item (Hemicontinuity) For any $t\in[0,T]$ and $x,y,z\in\mX$, the mapping
$$
[0,1]\ni\eps\mapsto [x,A(t, y+\eps z)]_\mX
$$
is continuous.

\item (Weak coercivity) There exist $q_1,q_2\geq 2$ and
$\lambda_1,\lambda_2, \lambda_3>0$ such that for all $x\in\mX$ and $t\in[0,T]$
\ce
[x,A(t,x)]_\mX\leq-\lambda_1\cdot\|x\|^{q_1}_{\mX_1}
-\lambda_2\cdot\|x\|^{q_2}_{\mX_2}+\lambda_3\cdot(\|x\|^2_\mH+1).
\de

\item (Weak monotonicity) There exist $\lambda_0, \lambda'_1, \lambda'_2\geq0$ such that for all
$x,y\in\mX$ and $t\in[0,T]$
\ce
[x-y,A(t,x)-A(t,y)]_\mX&\leq&-\lambda'_1\|x-y\|^{q_1}_{\mX_1}-\lambda'_2\|x-y\|^{q_2}_{\mX_2}\\
&&+\lambda_0\cdot\|x-y\|^2_\mH,
\de
where $q_1$ and $q_2$ are same as in ({\bf H}2).

\item (Boundedness) There exist $c_{A_1}, c_{A_2}>0$ such that for all $x\in\mX$ and $t\in[0,T]$
\ce
\|A_i(t,x)\|_{\mX^*_i}\leq c_{A_i}\cdot(\|x\|^{q_i-1}_{\mX_i}+1),\ \ i=1,2,
\de
where $q_1$ and $q_2$ are same as in ({\bf H}2).
\item There exists a $\beta_1>0$ such that for all $x,y\in\mH$ and $t\in[0,T]$
$$
\|B(t,x)-B(t,y)\|_{L_2(\mU_Q,\mH)}\leq \beta_1\|x-y\|_\mH
$$
and
$$
\|B(t,x)\|_{L_2(\mU_Q,\mH)}\leq \beta_1(1+\|x\|_\mH).
$$
\end{enumerate}

We take the polish space $\mS$ in Theorem \ref{Th2} as follows
\be
\mS:=\mC_T(\mH)\cap L^{q_1}(0,T;\mX_1)\cap L^{q_2}(0,T;\mX_2)\label{MS}
\ee
with the norm
$$
\|x\|_\mS:=\sup_{t\in[0,T]}\|x(t)\|_\mH+\sum_{i=1,2}
\left(\int^T_0\|x(t)\|^{q_i}_{\mX_i}\dif t\right)^{1/q_i}.
$$

Consider the following stochastic evolution equation:
\be
\label{SED}\left\{
\begin{array}{ll}
\dif X(t)=A(t,X(t))\dif t
+B(t,X(t))\dif W(t),\\
X(0)=x_0\in\mH.
\end{array}
\right.
\ee
By \cite{Kr-Ro, Zh1} and \cite{Roe-Zh}, we have the following existence of unique strong
solution to Eq.(\ref{SED}).
\bt\label{Th11}
Assume that {\bf (H1)}-{\bf (H5)}  hold. Then there exists a unique
measurable functional $\Phi$ from $\mC_T(\mU)$ to $\mS$ such that
$X(t,\omega)=\Phi(W_\cdot(\omega))(t)$ solves the following equation in $\mX^*$
$$
X(t)=x_0+\int^t_0A(s,X(s))\dif s+\int^t_0B(s,X(s))\dif W(s),
$$
where the It\^o stochastic integral is calculated  in Hilbert space $\mH$.
Moreover, for any $h\in\cA_N$, $X^h(t,\omega)=\Phi(W_\cdot(\omega)+h(\omega))(t)$
solves the following equation in $\mX^*$
$$
X(t)=x_0+\int^t_0A(s,X(s))\dif s+\int^t_0B(s,X(s))\dif W(s)+\int^t_0B(s,X(s))\dot{h}(s)\dif s.
$$
\et
\br
The second conclusion follows from the Girsanov theorem.
\er

\section{Laplace and Large Deviation Principle}

Consider the following small perturbation to stochastic evolution equation (\ref{SED}):
\be
\label{PSED}\left\{
\begin{array}{ll}
\dif X_\eps(t)=A(t,X_\eps(t))\dif t +\sqrt{\eps}B(t,X_\eps(t))\dif W(t),\ \eps\in(0,1),\\
X_\eps(0)=x_0\in\mH.
\end{array}
\right.
\ee
By Theorem \ref{Th11}, there exists a measurable mapping $\Phi_\eps: \mC_T(\mU)\to\mS$
such that
$$
X_\eps(t,\omega)=\Phi_\eps(W_\cdot(\omega))(t).
$$
We now fix a family of processes $\{h^\eps\}$ in $\cA_N$, and put
$$
X^\eps(t,\omega):=\Phi_\eps(W_\cdot(\omega)+\frac{h^{\eps}(\omega)}{\sqrt{\eps}})(t).
$$
It should be noticed that we have used a little confused notations $X_\eps$ and $X^\eps$,
but it is clearly different.
Note that $X^\eps(t)$ solves the following  stochastic evolution equation:
\be
\label{Eq1}\left\{
\begin{array}{ll}
\dif X^\eps(t)=A(t,X^\eps(t))\dif t
+\sqrt{\eps}B(t,X^\eps(t))\dif W(t)+B(t,X^\eps(t)){\dot h}^\eps(t)\dif t,\\
X^\eps(0)=x_0\in\mH.
\end{array}
\right.
\ee
Moreover, the following energy identity holds(cf. \cite{Kr-Ro}, also called It\^o's formula):
\be
\|X^\eps(t)\|_\mH^2&=&\|x_0\|_\mH^2+2\int^t_0[X^\eps(s),A(s,X^\eps(s))]_\mX\dif s+M^\eps(t)\no\\
&&+2\int^t_0\<X^\eps(s),B(s,X^\eps(s)){\dot h}^\eps(s)\>_\mH\dif s\no\\
&&+\eps\int^t_0\|B(s,X^\eps(s))\|^2_{L_2(\mU_Q,\mH)}\dif s,\label{Es}
\ee
where $t\mapsto M^\eps(t)$ is a real continuous martingale given by
$$
M^\eps(t):=2\sqrt{\eps}\int^t_0\<X^\eps(s), B(s,X^\eps(s))\dif W(s)\>_\mH.
$$
Note that the square variation process of $M^\eps(t)$ is given by
$$
<M^\eps>_t=4\eps\sum_{j}\int^t_0\<X^\eps(s), B(s,X^\eps(s))Q^{1/2}(e_j)\>_\mH^2\dif s,
$$
where $\{e_j\}$ is an orthogonormal basis of $\mU$.

\textsc{Convention}: The letter $C$ below with or without subscripts will
denote positive constants whose values may change in different occasions.

Our main task is to verify the above {\bf (Hypothesis)}.
We first prove some uniform estimates about $X^\eps(t)$.

\bl\label{Le1}
For any $p\geq 2$ and $T>0$, there exists a constant $C_{p,T}>0$ such that for all $\eps\in(0,1]$
$$
\mE\left(\sup_{t\in[0,T]}\|X^\eps(t)\|_\mH^{2p}\right)\leq C_{p,T}(\|x_0\|^{2p}_\mH+1),
$$
and for $i=1,2$
$$
\mE\left(\int^{T}_0\|X^\eps(s)\|^{q_i}_{\mX_i}\dif s\right)^p\leq C_{p,T}(\|x_0\|^{2p}_\mH+1).
$$
\el
\begin{proof}
By (\ref{Es}) and It\^o's formula, we find that
\ce
\|X^\eps(t)\|_\mH^{2p}&=&\|x_0\|_\mH^{2p}+2p\int^t_0\|X^\eps(s)\|_\mH^{2p-2}\cdot
[X^\eps(s),A(s,X^\eps(s))]_\mX\dif s\\
&&+p\int^t_0\|X^\eps(s)\|_\mH^{2p-2}\dif M^\eps(s)+\frac{p(p-1)}{2}\int^t_0
\|X^\eps(s)\|_\mH^{2(p-2)}\dif <M^\eps>_s\\
&&+2p\int^t_0\|X^\eps(s)\|_\mH^{2p-2}\<X^\eps(s),B(s,X^\eps(s)){\dot h}^\eps(s)\>_\mH\dif s\\
&&+\eps p\int^t_0\|X^\eps(s)\|_\mH^{2p-2}\|B(s,X^\eps(s))\|^2_{L_2(\mU_Q,\mH)}\dif s.
\de
By {\bf (H2)} and {\bf (H5)} we have
\ce
\|X^\eps(t)\|_\mH^{2p}&\leq&\|x_0\|_\mH^{2p}
+C\int^t_0\Big(\|X^\eps(s)\|_\mH^{2p}(\|{\dot h}^\eps(s)\|_{\mU_Q}+1)+1\Big)\dif s\\
&&+p\int^t_0\|X^\eps(s)\|_\mH^{2p-2}\dif M^\eps(s).
\de
Hence, by Gronwall's inequality and (\ref{Op2}) we get
$$
\|X^\eps(t)\|_\mH^{2p}\leq C_N\left(\|x_0\|_\mH^{2p}+1
+\int^t_0\left|\int^s_0\|X^\eps(r)\|_\mH^{2p-2}\dif M^\eps(r)\right|\dif s\right).
$$

Put
$$
f(t):=\mE\left(\sup_{r\in[0,t]}\|X^\eps(r)\|_\mH^{2p}\right).
$$
Then, by BDG's inequality and Young's inequality we obtain
\ce
f(t)&\leq& C_N\cdot(\|x_0\|_\mH^{2p}+1)
+C_N\cdot T\cdot\mE\left(\int^t_0\|X^\eps(r)\|_\mH^{4p-4}\dif <M^\eps>_r\right)^{1/2}\\
&\leq& C_N\cdot(\|x_0\|_\mH^{2p}+1)
+C_N\mE\left(\sup_{r\in[0,t]}\|X^\eps(r)\|_\mH^{2p}\cdot
\int^t_0(\|X^\eps(r)\|_\mH^{2p}+1)\dif r\right)^{1/2}\\
&\leq& C_N\cdot(\|x_0\|_\mH^{2p}+1)
+\frac{1}{2}f(t)+C_N\int^t_0(\mE\|X^\eps(r)\|_\mH^{2p}+1)\dif r.
\de
Therefore,
$$
f(t)\leq 2C_N\cdot(\|x_0\|_\mH^{2p}+1)+2C_N\int^t_0(f(r)+1)\dif r.
$$
By Gronwall's inequality again, we obtain the first estimate.

As for the second estimate, from (\ref{Es}) and {\bf (H2)}, {\bf (H5)}  we also have
\ce
\sum_{i=1,2}\int^{T}_0\|X^\eps(s)\|^{q_i}_{\mX_i}\dif s
\leq C\Big(\|x_0\|_\mH^2+|M^\eps(t)|+
\int^t_0\|X^\eps(s)\|^2_\mH(\|{\dot h}^\eps(s)\|_{\mU_Q}+1)\dif s\Big).
\de
Using the first estimate, we immediately get the desired second estimate.
\end{proof}

\bl\label{Le2}
For any $p\geq 2$, there exists a constant $C$ depending on $p,T,N$ and $x_0$ such that
for all $t,r\in[0,T]$ and $\eps\in(0,1)$
$$
\mE\|X^\eps(t)-X^\eps(r)\|^p_{\mX^*}\leq C|t-r|^{\frac{p}{q_1\vee q_2}}.
$$
\el
\begin{proof}
Note that the following equality holds in $\mX^*$
\ce
X^\eps(t)-X^\eps(r)&=&\int^t_rA(s,X^\eps(s))\dif s+
\sqrt{\eps}\int^t_rB(s,X^\eps(s))\dif W(s)\\
&&+\int^t_rB(s,X^\eps(s)){\dot h}^\eps(s)\dif s.
\de
Hence
$$
\mE\|X^\eps(t)-X^\eps(r)\|^p_{\mX^*}\leq 3^{p-1}(I_1+I_2+I_3),
$$
where
\ce
I_1&:=&\mE\left(\int^t_r\|A(s,X^\eps(s))\|_{\mX^*}\dif s\right)^p\\
I_2&:=&\mE\left\|\sqrt{\eps}\int^t_rB(s,X^\eps(s))\dif W(s)\right\|^p_{\mX^*}\\
I_3&:=&\mE\left\|\int^t_rB(s,X^\eps(s)){\dot h}^\eps(s)\dif s\right\|^p_{\mX^*}.
\de

For $I_1$, we have by {\bf (H4)} and H\"older's inequality
\ce
I_1&\leq&C\mE\left(\int^t_r(\|A_1(s,X^\eps(s))\|_{\mX^*_1}
+\|A_2(s,X^\eps(s))\|_{\mX^*_2})\dif s\right)^p\\
&\leq&C\mE\left(\int^t_r(\|X^\eps(s)\|^{q_1-1}_{\mX_1}+\|X^\eps(s)\|^{q_2-1}_{\mX_2}+1)\dif s\right)^p\\
&\leq&C(t-r)^p+C\sum_{i=1,2}\left[\mE\left(\int^t_r\|X^\eps(s)\|^{q_i}_{\mX_i}\dif s
\right)^{\frac{p(q_i-1)}{q_i}}(t-r)^{\frac{p}{q_i}}\right].
\de

For $I_2$, we have by BDG's inequality and  {\bf (H5)}
\ce
I_2&\leq&C\mE\left\|\sqrt{\eps}\int^t_rB(s,X^\eps(s))\dif W(s)\right\|^p_{\mH}\\
&\leq&C\mE\left(\int^t_r\|B(s,X^\eps(s))\|^2_{L_2(\mU_Q,\mH)}\dif s\right)^{p/2}\\
&\leq&C|t-r|^{p/2-1}\left(\int^t_r(\mE\|X^\eps(s))\|^p_{\mH}+1)\dif s\right).
\de

For $I_3$, we have by H\"older's inequality, (\ref{Op2})  and {\bf (H5)}
\be
I_3&\leq&C\mE\left\|\int^t_rB(s,X^\eps(s)){\dot h}^\eps(s)\dif s\right\|^p_{\mH}\no\\
&\leq&C\mE\left(\int^t_r\|B(s,X^\eps(s))\|_{L_2(\mU_Q,\mH)}\cdot
\|{\dot h}^\eps(s)\|_{\mU_Q}\dif s\right)^p\no\\
&\leq&C\mE\left(\int^t_r\|B(s,X^\eps(s))\|^2_{L_2(\mU_Q,\mH)}\dif s\right)^{p/2}
\left(\int^T_0\|{\dot h}^\eps(s)\|^2_{\mU_Q}\dif s\right)^{p/2}\no\\
&\leq&C|t-r|^{p/2-1}\left(\int^t_r(\mE\|X^\eps(s))\|^p_{\mH}+1)\dif s\right)\cdot N^{p/2}.\label{Es2}
\ee
The desired estimate now follows by combining the above estimates and Lemma \ref{Le1}.
\end{proof}

\bl\label{Le3}
Assume that for almost all $\omega$,
$h^\eps(\cdot,\omega)$ weakly converge to $h(\cdot,\omega)$ in $\mL_Q$,
and $X^\eps(\cdot,\omega)$ strongly converge to $X(\cdot,\omega)$ in $\mC_T(\mX^*)$.
Then $X(\cdot,\omega)$ solves the following equation
$$
X(t,\omega)=x_0+\int^t_0A(s,X(s,\omega))\dif s+\int^t_0B(s,X(s,\omega)){\dot h}(s,\omega)\dif s.
$$
Moreover, there exists a subsequence $\eps_k$ such that as $k\rightarrow\infty$,
\be
\mE\left(\sup_{t\in[0,T]}\|X^{\eps_k}(t)-X(t)\|^2_\mH\right)\rightarrow 0,\label{Lim1}
\ee
and if $\lambda_1',\lambda_2'>0$ in {\bf (H3)}, then for $i=1,2$
\be
\int^T_0\mE\|X^{\eps_k}(t)-X(t)\|^{q_i}_{\mX_i}\dif t\rightarrow 0.\label{Lim2}
\ee
\el
\begin{proof}
Set for $i=1,2$
$$
\mK_{1,i}:=L^{q_i/(q_i-1)}([0,T]\times\Omega, \cM,\dif t\times\dif P;\mX^*_i)
$$
and
$$
\mK_{2,i}:=L^{q_i}([0,T]\times\Omega, \cM,\dif t\times\dif P;\mX_i),
$$
where $\cM$ denotes the progressively $\sigma$-algebra associated with $\cF_t$.
Then $\mK_{1,i}$ and $\mK_{2,i}$ are reflexive and separable Banach spaces.

We have by Lemma \ref{Le1}
\be
\sup_{\eps\in(0,1]}\mE\|X^\eps(T)\|_\mH^{2}+\sup_{\eps\in(0,1]}
\sum_{i=1,2}\|X^\eps\|^{q_i}_{\mK_{2,i}}
<+\infty\label{PP1}
\ee
and
\be
\sup_{\eps\in(0,1]}\mE\left(\sup_{t\in[0,T]}\|X^\eps(t)\|_\mH^4\right)<+\infty.\label{PP6}
\ee
Hence, by the strong convergence of $X^\eps(\cdot,\omega)$ to
$X(\cdot,\omega)$ in $\mC_T(\mX^*)$ we have
\be
\mE\|X(T)\|_\mH^{2}&\leq&\varliminf_{\eps\downarrow 0}\mE\|X^\eps(T)\|_\mH^{2}
\leq C_{p,T}(\|x_0\|^{2}_\mH+1)\label{Lim9}\\
\int^{T}_0\mE\|X(s)\|^{q_i}_{\mX_i}\dif s&\leq&\varliminf_{\eps\downarrow 0}
\int^{T}_0\mE\|X^\eps(s)\|^{q_i}_{\mX_i}\dif s
\leq C_{p,T}(\|x_0\|^{2}_\mH+1),
\ee
as well as by (\ref{PP6})
$$
\lim_{\eps\downarrow 0}\mE\left(\sup_{t\in[0,T]}\|X^\eps(s)-X(s)\|^2_{\mX^*}\right)=0.
$$
Thus, by (\ref{PP4}) we have as $\eps\downarrow 0$
\be
&&\mE\left(\int^T_0\|X^\eps(s)-X(s)\|^2_\mH\dif s\right)
=\mE\left(\int^T_0[X^\eps(s)-X(s),X^\eps(s)-X(s)]_\mX\dif s\right)\no\\
&\leq&\int^T_0\mE\Big(\|X^\eps(s)-X(s)\|_\mX\cdot \|X^\eps(s)-X(s)\|_{\mX^*}\Big)\dif s\no\\
&\leq&\left(\int^T_0\mE\|X^\eps(s)-X(s)\|^2_\mX\dif s\cdot
\int^T_0\mE\|X^\eps(s)-X(s)\|^2_{\mX^*}\dif s\right)^{1/2}\rightarrow 0.\label{PP2}
\ee

Notice also that by {\bf (H4)} and (\ref{PP1})
$$
\sup_{\eps\in(0,1]}\|A_i(\cdot,X^\eps(\cdot))\|_{\mK_1,i}<+\infty,\ \ i=1,2.
$$
By this and (\ref{PP1}) and the weak compactness of $\mK_{1,i}$ and $\mK_{2,i}$, $i=1,2$,
there exist a subsequence $\eps_k$(still denoted by $\eps$ for simplicity) and
$Y_i\in \mK_{1,i}, i=1,2$, $\bar X\in \cap_{i=1,2}\mK_{2,i}$, $X_T\in L^2(\Omega)$ such that
\be
X^\eps&\rightarrow& \bar X
\mbox{ weakly in $\mK_{2,i}$}, i=1,2,\label{Lim4}\\
X^\eps(T)&\rightarrow& X_T\mbox{ weakly in $L^2(\Omega)$},
\ee
and
\be
Y^\eps_i:=A_i(\cdot,X^\eps(\cdot))\rightarrow Y_i
\mbox{ weakly in $\mK_{1,i}$}, i=1,2.\label{Lim5}
\ee

Put $Y=Y_1+Y_2$ and define
$$
\tilde X(t):=x_0+\int^t_0Y(s)\dif s+\int^t_0B(s,X(s)){\dot h}(s)\dif s.
$$
Note that
\ce
X^\eps(t)&=&x_0+\int^t_0A(s,X^\eps(s))\dif s+
\sqrt{\eps}\int^t_0B(s,X^\eps(s))\dif W^\eps(s)\\
&&+\int^t_0B(s,X^\eps(s))\dot{ {h}}^\eps(s)\dif s.
\de
By taking weak limits and (\ref{PP2}), it is not hard to see that
(see also the proof of (\ref{Lim3}) below)
$$
\tilde X(t,\omega)=X(t,\omega)=\bar X(t,\omega)
 \mbox{ for $\dif t\times \dif P$-almost all $(t,\omega)$}
$$
and
\be
\tilde X(T)=X(T)=X_T.\label{Op3}
\ee

In the following we use the unified notation $X$, and only need to prove
by the usual monotonicity argument that
\be
Y(s,\omega)=A(s,X(s,\omega))\mbox{ for $\dif t\times \dif P$-almost all $(t,\omega)$}.\label{Op1}
\ee
Without loss of generality, we assume that $\lambda_0=0$ in {\bf (H3)}(cf. \cite{Roe, Zh1}).
It is clear that in (\ref{Es})
\be
\mE M^\eps(T)=0, \ \label{Lim8}
\ee
and
\be
\lim_{\eps\downarrow 0} \eps\mE\int^T_0\|B(s,X^\eps(s))\|^2_{L_2(\mU_Q,\mH)}\dif s=0.
\ee
Let us prove the following limit:
\be
\lim_{\eps\downarrow 0}\mE\left|\int^T_0\Big(\<X^\eps(s),B(s,X^\eps(s)){\dot h}^\eps(s)\>_\mH
-\<X(s),B(s,X(s)){\dot h}(s)\>_\mH\Big)\dif s\right|=0.\label{Lim3}
\ee
Since for almost all $\omega$,
$h^\eps(\cdot,\omega)$ weakly converges to $h(\cdot,\omega)$ in $\mL_Q$, by the dominated convergence
theorem we have
$$
\lim_{\eps\downarrow 0}\mE\left|\int^T_0\<X(s),B(s,X(s))({\dot h}^\eps(s)
-{\dot h}(s))\>_\mH\dif s\right|=0.
$$
By (\ref{Op2}), Lemma \ref{Le1} and (\ref{PP2}) we also have
\be
&&\mE\left|\int^T_0\<X^\eps(s)-X(s),B(s,X^\eps(s)){\dot h}^\eps(s)\>_\mH\dif s\right|\no\\
&\leq&\mE\int^T_0\|X^\eps(s)-X(s)\|_{\mH}\cdot
(\|X^\eps(s)\|_\mH+1)\cdot\|{\dot h}^\eps(s)\|_{\mU_Q}\dif s\no\\
&\leq&C_N\mE\left(\int^T_0\|X^\eps(s)-X(s)\|_{\mH}^2\cdot(\|X^\eps(s)\|_\mH+1)^2\dif s\right)^{1/2}\no\\
&\leq&C_N\mE\left(\Big(\sup_{s\in[0,T]}\|X^\eps(s)\|^2_\mH+1\Big)\cdot
\int^T_0\|X^\eps(s)-X(s)\|_{\mH}^2\cdot\dif s\right)^{1/2}\no\\
&\leq&C_N\left(\int^T_0\mE\|X^\eps(s)-X(s)\|_{\mH}^2\dif s\right)^{1/2}
\rightarrow 0,\label{PP3}
\ee
and
\ce
&&\mE\left|\int^T_0\<X(s),(B(s,X^\eps(s))-B(s,X(s))){\dot h}^\eps(s)\>_\mH\dif s\right|\\
&\leq&\mE\int^T_0\|X(s)\|_\mH\cdot \|X^\eps(s)-X(s)\|_{\mH}
\cdot\|{\dot h}^\eps(s)\|_{\mU_Q}\dif s\\
&\leq&C_N\mE\left(\int^T_0\|X(s)\|_\mH^2\cdot \|X^\eps(s)-X(s)\|_{\mH}^2
\dif s\right)^{1/2}\rightarrow 0.
\de
The limit (\ref{Lim3}) now follows.

Notice that for any $\Phi\in\mK_{2,1}\cap\mK_{2,2}$
\be
\mE\int^T_0[ X^{\eps}(s),A(s, X^{\eps}(s))]_\mX\dif s
&\leq&\mE\int^T_0[\Phi(s),A(s,  X^{\eps}(s))-A(s,\Phi(s))]_\mX\dif s\no\\
&&+\mE\int^T_0[ X^{\eps}(s),A(s,\Phi(s))]_\mX\dif s\ \ (\because \lambda_0=0)\no\\
&\rightarrow&\mE\int^T_0[\Phi(s),Y(s)-A(s,\Phi(s))]_\mX\dif s\no\\
&&+\mE\int^T_0[ X(s),A(s,\Phi(s))]_\mX\dif s,\label{Lim6}
\ee
as $\eps\downarrow 0$. The limits are due to (\ref{Lim4}) and  (\ref{Lim5}).

Combining (\ref{Es}), (\ref{Lim9}), (\ref{Lim8})-(\ref{Lim6}) yields that
\ce
&&\mE\|X_T\|_\mH^{2}\leq\varliminf_{\eps\downarrow 0}\mE\|X^\eps(T)\|_\mH^{2}\\
&\leq& \|x_0\|_\mH^{2}+2\mE\int^T_0[\Phi(s),Y(s)-A(s,\Phi(s))]_\mX\dif s\\
&&+2\mE\int^T_0[ X(s),A(s,\Phi(s))]_\mX\dif s\\
&&+2\mE\int^T_0\<X(s),B(s,X(s)){\dot h}(s)\>_\mH\dif s.
\de
On the other hand, by the energy equality(see (\ref{Es})) we have
\ce
\|\tilde X(T)\|_\mH^{2}=\|x_0\|_\mH^{2}+2\int^T_0[X(s),Y(s)]_\mX\dif s
+2\int^T_0\<X(s),B(s,X(s)){\dot h}(s)\>_\mH\dif s.
\de
So by (\ref{Op3})
\ce
\mE\int^T_0[X(s)-\Phi(s),Y(s)-A(s,\Phi(s))]_\mX\dif s\leq 0,
\de
which then yields (\ref{Op1}) by {\bf (H1)} and \cite[Lemma 2.5]{Zh1}(see also \cite{Roe}).

Lastly, let us prove the limits (\ref{Lim1}) and (\ref{Lim2}). By It\^o's formula, we have
\ce
\|X^{\eps}(t)-X(t)\|_\mH^2=I^{\eps}_1(t)+I^{\eps}_2(t)+I^{\eps}_3(t)+I^{\eps}_4(t).
\de
where
\ce
I^\eps_1(t)&:=&2\int^t_0[X^{\eps}(s)-X(s),A(s,X^{\eps}(s))-A(s,X(s))]_\mX\dif s\\
I^\eps_2(t)&:=&2\int^t_0\<X^{\eps}(s)-X(s),B(s,X^{\eps}(s)){\dot h}^{\eps}(s)\>_\mH\dif s\no\\
I^\eps_3(t)&:=&-2\int^t_0\<X^{\eps}(s)-X(s),B(s,X(s)){\dot h}(s)\>_\mH\dif s\no\\
I^\eps_4(t)&:=&2\sqrt{\eps}\int^t_0\<X^{\eps}(s)-X(s),B(s,X^{\eps}(s))\dif W(s)\>_\mH\no\\
I^\eps_5(t)&:=&\eps\int^t_0\|B(s,X^{\eps}(s))\|^2_{L_2(\mU_Q,\mH)}\dif s.
\de

By BDG's inequality and Lemma \ref{Le1}, we obviously have
\ce
\lim_{k\rightarrow\infty}\mE\left(\sup_{t\in[0,T]}(|I^\eps_4(t)|+|I^\eps_5(t)|)\right)=0.
\de

For $I^\eps_2$, as in the proof of (\ref{PP3}) we have
\ce
\mE\left(\sup_{t\in[0,T]}|I^\eps_2(t)|\right)&\leq&
C\mE\int^T_0\|X^{\eps}(s)-X(s)\|_\mH\cdot
(\|X^{\eps}(s)\|_\mH+1)\cdot\|{\dot h}^{\eps}(s)\|_{\mU_Q}\dif s\\
&\leq&C_N\mE\left(\int^T_0\|X^{\eps}(s)-X(s)\|^2_\mH\cdot
(\|X^{\eps}(s)\|_\mH+1)^2\dif s\right)^{1/2}\\
&\rightarrow& 0,\ \ \ \mbox{  as $\eps\rightarrow 0$}.
\de
Similarly
$$
\lim_{k\rightarrow\infty}\mE\left(\sup_{t\in[0,T]}|I^\eps_3(t)|\right)=0.
$$
Assume $\lambda_1',\lambda_2'>0$, then
\be
I^\eps_1(t)\leq -\sum_{i=1,2}\lambda_i'\int^t_0\|X^{\eps}(s)-X(s)\|^{q_i}_{\mX_i}\dif s+
\lambda_0\int^t_0\|X^{\eps}(s)-X(s)\|^2_{\mH}\dif s.\label{Es5}
\ee
If we put
$$
f(t):=\varlimsup_{\eps\rightarrow\infty}\mE\left(\sup_{s\in[0,t]}\|X^{\eps}(s)-X(s)\|_\mH^2\right),
$$
then
$$
f(t)\leq \lambda_0\int^t_0f(s)=0.
$$
So
$$
f(T)=0.
$$
The limits (\ref{Lim1}) and (\ref{Lim2}) is straightforward by noting (\ref{Es5}).
\end{proof}

We may prove the following main lemma.
\bl\label{Le4}
There exists a probability space
$(\tilde\Omega,\tilde\cF,\tilde P)$ and a sequence
(still indexed by $\eps$ for simplicity) $\{(\tilde h^\eps, \tilde X^\eps,\tilde W^\eps)\}$
and $(h, X^h, \tilde W)$ defined on this probability
space and taking values in $D_N\times \mC_T(\mX^*)\times \mC_T(\mU)$ such that

(a) $(\tilde h^\eps, \tilde X^\eps,\tilde W^\eps)$ has the same law as
$(h^\eps, X^\eps,W)$ for each $\eps$;

(b) $(\tilde h^\eps, \tilde X^\eps,\tilde W^\eps)
\rightarrow (h, X^h,\tilde W)$ in $D_N \times \mC_T(\mX^*)\times \mC_T(\mU)$,
$\tilde P$-a.s. as $\eps\rightarrow 0$;

(c) $(h,X^h)$ uniquely solves the following equation:
\be
X^h(t)=x_0+\int^t_0A(s,X^h(s))\dif s+\int^t_0B(s,X^h(s)){\dot h}(s)\dif s.\label{Es4}
\ee

Moreover, there exists a subsequence $\eps_k$ such that as $k\rightarrow\infty$,
\be
\mE^{\tilde P}\left(\sup_{t\in[0,T]}\|\tilde X^{\eps_k}(t)
-X^h(t)\|^2_\mH\right)\rightarrow 0,\label{Lim10}
\ee
and if $\lambda_1',\lambda_2'>0$ in {\bf (H3)}, then for $i=1,2$
\be
\int^T_0\mE^{\tilde P}\|\tilde X^{\eps_k}(t)-X^h(t)\|^{q_i}_{\mX_i}\dif t\rightarrow 0.\label{Lim20}
\ee
\el
\begin{proof}
By Lemma \ref{Le2} and \cite[Corollary 14.9]{Ka}, the laws of  $(h^\eps, X^\eps, W)$
in $D_N\times \mC_T(\mX^*)\times \mC_T(\mU)$ is tight. By Skorohod's embedding theorem,
the conclusions (a) and (b) hold.
Note that $\tilde X^\eps(0)=x_0$ $\tilde P$-a.s. and
\ce
\tilde X^\eps(t)&=&x_0+\int^t_0A(s,\tilde X^\eps(s))\dif s+
\sqrt{\eps}\int^t_0B(s,\tilde X^\eps(s))\dif \tilde W^\eps(s)\\
&&+\int^t_0B(s,\tilde X^\eps(s))\dot{ \tilde{h}}^\eps(s)\dif s.
\de
The other conclusions follows from Lemma \ref{Le3}.
\end{proof}

From this lemma, one sees that {\bf (Hypothesis)} holds. Thus, by Theorem \ref{Th2} we obtain
\bt\label{Th1}
Assume {\bf (H1)}-{\bf (H5)} hold, and $\lambda_1',\lambda_2'>0$ in {\bf (H3)}.
Then for all real bounded continuous functions $g$ on $\mS$
$$
\lim_{\eps\rightarrow 0}\eps\log\mE\left(\exp\left[-\frac{g(X_\eps)}{\eps}\right]\right)
=-\inf_{f\in\mS}\{g(f)+I(f)\},
$$
where $I(f)$ is defined by
\be
I(f):=\frac{1}{2}\inf_{\{h\in\mL_Q:~f=X^h\}}\|h\|^2_{\mL_Q},\label{Rate}
\ee
and $X^h$ solves (\ref{Es4}).
\et
\br
If $\lambda_1',\lambda_2'=0$ in {\bf (H3)}, then the conclusion still holds
if $\mS$ is replaced by $\mC_T(\mH)$.
\er

In order to show the large deviation principle, we need to prove that $I(f)$
is a good rate function. For this aim, we need an extra assumption:
$$
\mX\hookrightarrow\mH \mbox{ compactly}.
$$

\bl\label{Th2}
In addition to {\bf (H1)}-{\bf (H5)} and $\lambda_1',\lambda_2'>0$, we also assume that
$\mX$ is compactly embedded in $\mH$. Then $I(f)$ is a good rate function, i.e.,
for any $a>0$, $\{f\in\mS: I(f)\leq a\}$ is compact.
\el
\begin{proof}
It suffices to prove that if $h_n\in D_N$ weakly converge to $h$ in $\mL_Q$, then
there exists a subsequence $n_k$ (still denoted by $n$) such that
\be
\lim_{n\rightarrow\infty}\|X^{h_n}-X^h\|_\mS=0.\label{Lim0}
\ee
In fact, assume that $I(f_n)\leq a$. By the definition of $I(f_n)$, there exists a
sequence $h_n\in\mL_Q$ such that $X^{h_n}=f_n$ and
$$
\frac{1}{2}\|h_n\|_{\mL_Q}^2\leq a+\frac{1}{n}.
$$
By the weak compactness of $D_{2a+1}$, there exists a subsequence $n_k$(still denoted by $n$)
and $h\in\mL_Q$ such that $h_n$ weakly converge to $h$ and
$$
\|h\|_{\mL_Q}^2\leq\varliminf_{n\rightarrow\infty}\|h_n\|_{\mL_Q}^2\leq 2a.
$$
Thus, by (\ref{Lim0}) we get the desired compactness.

We now prove (\ref{Lim0}). As in the proofs of Lemma \ref{Le1} and Lemma \ref{Le2}, we may prove
$$
\|X^{h_n}\|_\mS=\sup_{t\in[0,T]}\|X^{h_n}(t)\|_\mH+\sum_{i=1,2}
\left(\int^T_0\|X^{h_n}(t)\|^{q_i}_{\mX_i}\dif t\right)^{1/q_i}\leq C
$$
and
$$
\|X^{h_n}(t)-X^{h_n}(r)\|_{\mX^*}\leq C|t-r|^{\frac{1}{q_1\vee q_2}},
$$
where $C$ is independent of $n$.

Since $\mX\hookrightarrow\mH$ is compact, by \cite[Theorem 2.1]{Fl-Ga}  there exists a subsequence $n_k$
(still denoted by $n$) and an $X\in L^2(0,T;\mH)$ such that
$$
\int^T_0\|X^{h_n}(t)-X(t)\|_\mH^2\dif t=0.
$$
Basing on this convergence, as in the proof of Lemma \ref{Le3}, we in fact have $X=X^h$
and the desired limit (\ref{Lim0}) hold.
\end{proof}

Using Theorem \ref{Th1} and Lemma \ref{Th2}, we obtain the following large deviation
principle.
\bt\label{Main}
Assume {\bf (H1)}-{\bf (H5)} hold, and $\mX$ is compactly embedded in $\mH$,
$\lambda_1',\lambda_2'>0$ in {\bf (H3)}.
Let the law of $X_\eps$ in $\mS$ be denoted by $\nu_\eps$.
Then for any $A\in\sB(\mS)$
\ce
-\inf_{f\in A^o}I(f)\leq\liminf_{\eps\rightarrow 0}\eps\log\nu_\eps(A)
\leq\limsup_{\eps\rightarrow 0}\eps\log\nu_\eps(A)\leq -\inf_{f\in \bar A}I(f),
\de
where the closure and the interior are taken in $\mS$, and $I(f)$ is a good rate function defined by
(\ref{Rate}).
\et

\section{Applications}

\subsection{SDE with monotone drift}

We consider the following small perturbation of
stochastic ordinary differential equation with monotone
drift:
\ce
\dif X_\eps(t)=b(t,X_\eps(t))\dif t+\sqrt{\eps}\sigma(t,X_\eps(t))\dif W(t),\ \ X(0)=x_0\in\mR^d,
\de
where $W$ is an $m$-dimensional Brownian motion, $\sigma$ and $b$ satisfy that
\begin{enumerate}[({\bf H$\sigma$})]
\item There exists a constant $C_\sigma>0$ such that for all $t\in[0,T]$ and $x\in\mR^d$
$$
\|\sigma(t,x)-\sigma(t,y)\|_{\mR^{d\times m}}\leq C_\sigma\|x-y\|_{\mR^d}.
$$
\end{enumerate}
\begin{enumerate}[({\bf Hb})]
\item There exists a constant $C_b>0$ such that for all $t\in[0,T]$ and $x\in\mR^d$
$$
\<x-y, b(t,x)-b(t,y)\>_{\mR^d}\leq C_b\|x-y\|_{\mR^d}^2.
$$
\end{enumerate}

Let $\nu_\eps$ be the law of $X_\eps(t)$ in the continuous functions space $\mC_T(\mR^d)$.
Then the conclusion of Theorem \ref{Main} holds.

The following two examples can also be found in \cite{Zh2}.
\subsection{Stochastic reaction diffusion equation}

Let $\cO$ be an open and bounded set in Euclidean space $\mR^d$, where the boundary $\p\cO$ of $\cO$
is assumed to be smooth. For $q\geq 2$, let $W^{1,q}_0(\cO)$ and $W^{-1,\frac{q}{q-1}}(\cO)$ be the usual
Sobolev spaces(cf. \cite{Ad}).

Suppose that for each integer $j=1,\cdots,d$, we are given a  function $a_j:\cO\times\mR\mapsto\mR$ such that
\be
&&r\mapsto a_j(\xi,r)\mbox{ is continuous and non-decreasing for each $\xi\in\cO$} ,\label{PI1}\\
&&a_j(\xi,r)\xi\geq C_1|r|^{q_1}-C_2,\ \ (\xi,r)\in\cO\times\mR,\label{PI2}\\
&&|a_j(\xi,r)|\leq C_3(|r|^{q_1-1}+1),\ \ (\xi,r)\in\cO\times\mR,\label{PI3}
\ee
where $q_1\geq 2$ and $C_1,C_2,C_3>0$.

Let $b$ be another continuous function  satisfying (\ref{PI1})-(\ref{PI3}) and
with a different constant $q_2\geq 2$. Let $l^2$
be the usual Hilbert space of square summable real number sequences.
Let $\sigma(\xi,r):\cO\times\mR\mapsto l^2$
satisfy that $\sigma(\cdot,0)\in L^2(\cO;l^2)$ and for some $c_1>0$
$$
\|\sigma(\xi,r)-\sigma(\xi,r')\|_{l^2}\leq c_1\cdot|r-r'|,\quad \xi\in\cO, r,r'\in\mR,
$$

We consider the following small perturbation of
stochastic reaction diffusion equation
\be
\left\{
\begin{array}{ll}
\dif X_\eps(t,\xi)=\Big[\sum^d_{i=1}\p_i a_i(\xi,\p_i X_\eps(t,\xi))-b(\xi,X_\eps(t,\xi))\Big]\dif t\\
\quad\quad\quad\quad\ \ +\sqrt{\eps}\sum_{j=1}^\infty \sigma_j(\xi,X_\eps(t,\xi))\dif W_j(t),\\
X_\eps(t,\xi)=0,\ \ \forall \xi\in\p\cO,\\
X_\eps(0,\xi)=x(\xi)\in L^2(\cO),
\end{array}
\right.\label{Eq8}
\ee
where $W_j(t)=\<W(t),\ell_j\>_{\mU_Q}$ and $\{\ell_j,j\in\mN\}$
is an orthogonal basis of $\mU_\mQ$.

Set
$$
\mX_1:=W^{1,q_1}_0(\cO), \ \
\mX_2:=L^{q_2}(\cO),\ \
\mH:=L^2(\cO)
$$
and
$$
\mX_1^*:=W^{-1,\frac{q_1}{q_1-1}}(\cO),\ \ \mX_2^*:=L^{\frac{q_2}{q_2-1}}(\cO).
$$
Then
\ce
\mX:=\mX_1\cap\mX_2\subset\mH\subset(\mX_1^*+\mX_2^*)\subset\mX^*
\de
forms an evolutional triple.

Now, define for $u,v\in \mX_1$
$$
[A_1(u),v]_{\mX_1}:=-\sum^d_{i=1}\int_\cO a_i(\xi,\p_i u(\xi))\cdot \p_i v(\xi)\dif \xi
$$
and for $u,v\in\mX_2$
$$
[A_2(u),v]_{\mX_2}:=-\int_\cO b(\xi,u(\xi))\cdot v(\xi)\dif \xi.
$$
Clearly, for each $u\in\mX_1$, $[A_1(u),\cdot]_{\mX_1}\in\mX_1^*$ and for each
$u\in\mX_2$, $[A_2(u),\cdot]_{\mX_2}\in\mX_2^*$.
Thus,
$$
A_1: \mX_1\mapsto\mX_1^*,\ \ A_2:\mX_2\mapsto\mX_2^*,
$$
and it is easy to verify that
$A:=A_1+A_2$ satisfies {\bf (H1)}-{\bf (H4)}.

Moreover, if we define for $x\in\mH=L^2(\cO)$
$$
B(t,x):=\sum_{j=1}^\infty \sigma_j(\cdot,x(\cdot))\ell_j
$$
then for any $x,y\in L^2(\cO)$
\ce
\|B(t,x)-B(t,y)\|_{L_2(\mU_Q;\mH)}^2&=&\sum_{j=1}^\infty \|\sigma_j(\cdot,x(\cdot))
-\sigma_j(\cdot,y(\cdot))\|^2_{L^2(\cO)}\\
&\leq&C\|x(\cdot)-y(\cdot)\|^2_{L^2(\cO)}.
\de
Thus, {\bf (H5)} holds.

Let $\nu_\eps$ be the law of $X_\eps(t)$ in  $\mS$, where $\mS$
is defined by (\ref{MS}). Then the conclusion of Theorem \ref{Main} holds.

\subsection{Stochastic Porous Medium Equation}

As in the previous subsection, we consider the bounded domain $\cO$ in $\mR^d$
with smooth boundary.

For $p\geq 2$, set
$$
\mX:=L^p(\cO), \ \ \mH:=W^{-1,2}(\cO),\ \ \mX^*:=L^{p/(p-1)}(\cO).
$$
The inner product in $\mH$ is given by
$$
\<x,y\>_\mH:=\int_{\cO}(-\Delta)^{-1/2}x(\xi)\cdot(-\Delta)^{-1/2}y(\xi)~\dif \xi,\ \
x,y\in\mH=W^{-1,2}(\cO).
$$
Note that $-\Delta$ establishes an isomorphism between $W^{1,2}_0(\cO)$
and $W^{-1,2}(\cO)$. We shall identify $W^{1,2}_0(\cO)$ with the dual space $\mH^*$ of $\mH$,
and hence $\mH^*=W^{1,2}_0(\cO)\subset L^{p/(p-1)}(\cO)$. Thus, we have the evolution triple
$$
\mX\subset\mH\simeq\mH^*\subset \mX^*
$$
where $\simeq$ is understood through $-\Delta$.

Let $\phi_p(r):=r|r|^{{p-2}/2}$, and define for $x\in\mX=L^p(\cO)$
$$
A(x):=\Delta\phi_p(x).
$$
Then $A(x)\in\mX^*$ and {\bf (H1)}-{\bf (H4)} hold(cf. \cite{Sh, Roe}).

Let $B_1,\cdots, B_n\in L_2(\mU_Q,\mH)$, and define
\be
B(t,x):=\sum_{k=1}^ng_k([e_{n_1},x]_\mH,\cdots,[e_{n_k},x]_\mH)B_k, \  \ e_{n_j}\in\mH,\label{Cyl}
\ee
where $g_k$ are Lipschitz continuous functions on $\mR^{n_k}$.
Then such $B(t,x)$ satisfies {\bf (H5)}.
It should be noticed that if $\sigma\in C^\infty_b(\mR)$ is not linear, the mapping
$x\mapsto\sigma(x)$ is in general not Lipschitz from $W^{-1,2}(\cO)$ to $W^{-1,2}(\cO)$.

Consider the following small perturbation of
stochastic porous medium equation
\be
\left\{
\begin{array}{ll}
\dif X_\eps(t)=\Delta(\phi_p(X_\eps(t)))\dif t
+\sqrt{\eps}B(t,X_\eps(t))\dif W(t),\\
X_\eps(t,\xi)=0,\ \ \forall \xi\in\p\cO,\\
X_\eps(0,\xi)=x(\xi)\in W^{-1,2}(\cO).
\end{array}
\right.\label{Eq9}
\ee
Let $\nu_\eps$ be the law of $X_\eps(t)$ in  $\mC_T(\mH)\cap L^p(0,T;\mX)$.
Then the conclusion of Theorem \ref{Main} holds.

\vspace{5mm}

{\bf Acknowledgements:}

The second named author would like to thank Professor Benjamin Goldys for
providing him an excellent environment to work in the University of New South Wales.
His work is supported by ARC Discovery grant DP0663153 of Australia.
The authors also thanks Dr. Wei Liu for sending us his preprint \cite{Li}.


\begin{thebibliography}{999}

\bibitem{Ad}R.A. Adams: {\it Sobolev space}. Academic Press, 1975.

\bibitem{Az}R. Azencott: Grandes d\'eviations et applications.
{\it Ecole d'Edt\'e de Probabilit\'es de
Saint-Flour VIII, 1978, Lect. Notes. in Math.}, 779, 1-176. Springer, New York, 1980.

\bibitem{bd}M. Bou\'e and P. Dupuis: A variational representation for certain
functionals of Brownian motion. {\it Ann. of Prob.}, 1998, Vol.
26, No.4, 1641-1659.

\bibitem{bde}M. Bou\'e, P. Dupuis and R.S. Ellis: Large deviations for small noise
difusions with discontinuous statistics. {\it Prob. Theory Relat.Fields.}, 116,125-149(2000).

\bibitem{bd0}A. Budhiraja and P. Dupuis: A variational representation
for positive functionals of infinite dimensional Brownian motion.
{\it Probab. Math. Statist.} 20 (2000), no. 1,
Acta Univ. Wratislav. No. 2246, 39--61.

\bibitem{Bu-Du-Ma}A. Budhiraja, P. Dupuis and V. Maroulas:
Large deviations for infinite dimensional stochastic dynamical systems. to appear in
{\it Ann. of Prob.}.

\bibitem{Ce-Ro}S. Cerrai and M. R\"ockner: Large deviations for stochastic reaction
diffusion systems with multiplcative noise and non-Lipschitz reaction term.
{\it The Annals of Probability}, 32:1100-1139, 1996.

\bibitem{Pa1}E. Pardoux: Stochastic partial differential equations
and filtering of diffusion processes.
{\it Stochastic.} 1979, 127-167.

\bibitem{DR1}G. Da Prato and M. R\"ockner: Weak solutions to stochastic
porous media equations, {\it J. Evolution Equ.} 4(2004), 249--271.

\bibitem{DaZa}G. Da Prato and  J. Zabczyk: {\it Stochastic equations in infinite dimensions}.
Cambridge: Cambridge University Press, 1992.

\bibitem{de}P. Dupuis and R.S. Ellis: {\it A Weak Convergence Approach to the
Theory of Large Deviations}. Wiley, New-York, 1997.

\bibitem{Fe-Ku}J. Feng and T.G. Kurtz: {\it Large deviations for stochastic processes}.
Math. Serveys and Mono., Vol. 131, AMS(2006).

\bibitem{Fl-Ga}F. Flandoli and D. Gatarek: Martingale and stationary solutions for stochastic Navier-Stokes equations.
{\it Probab. Theory Related Fields}, 102, 367-391(1995).

\bibitem{Fr-We}M.I. Freidlin and A.D. Wentzell: On small random
perturbations of dynamical system, {\it Russian Math. Surveys} 25
(1970), 1-55.

\bibitem{Go-Mi}I. Gy\"ongy and A. Millet: On Discretization Schemes for Stochastic Evolution Equations.
{\it Potential Analysis}, (2005)23:99-134.

\bibitem{Ka-Xi}G. Kallianpur, J. Xiong: Large deviations for a class of stochastic
partial differential equations. {\it The Annals of Probability}, 24(1):320-345, 1996.

\bibitem{Ka}O. Kallenberg: {\it Foundations of Mordern Probability}. Springer-Verlag, Berlin, 1997.

\bibitem{Kr-Ro}N.V. Krylov and B.L. Rozovskii: Stochastic evolution equations.
{\it J. Soviet Math.}(Russian), 1979, pp. 71-147, Transl. 16(1981), 1233-1277.

\bibitem{Li}W. Liu: Large deviations for stochastic evolution equations
with small multiplcative noise. Preprint.

\bibitem{Pe}S. Peszat: Large deviation estimates for stochastic evolution equations.
{\it Prob. Th. Rel. Fields}, 98:113-136, 1994.

\bibitem{Roe}C. Pr\'ev\^ot  and M. R\"ockner: A concise course on stochastic partial
differential equations. Lecture Notes in Mathematics, 1905. Springer, Berlin, 2007. vi+144 pp

\bibitem{Re-Ro-Wa}J. Ren, M. R\"ockner and F. Wang:
Stochastic Generalized Porous Media and Fast Diffusion Equations. Preprint.


\bibitem{Rz1}J. Ren and X. Zhang: Schilder theorem for the Brownian motion on the
diffeomorphism group of the circle. {\it J. Func. Anal.}, Vol. 224, I. 1, 107-133(2005).

\bibitem{Rz2}J. Ren and X. Zhang: Freidlin-Wentzell's large deviations
for homeomorphism flows of non-Lipschitz SDEs. {\it Bull. Sci. Math. 2 Serie},
Vol 129/8 pp 643-655(2005).


\bibitem{Roe-Zh}M. R\"ockner, B. Schmuland and X. Zhang: Yamada-Watanabe Theorem for
Stochastic Evolution Equations in Infinite Dimensions. Preprint.

\bibitem{Ro-Wa-Wu}M. R\"ockner, F.Y. Wang and L. Wu: Large deviations for stochastic
generalized porous media equations. {\it Stoch. Proc. and their Appl.}, 116(2006)1677-1689.

\bibitem{Ro}B.L. Rozovskii: {\it Stochastic evolution systems. Linear theory and applications to
nonlinear filtering.} Mathematics and its Applications (Soviet Series), 35, Kluwer Academic Publishers, 1990.

\bibitem{Sh}R.E. Showalter: {\it Monotone Operators in Banach Space and Nonlinear Partial Differential
Equations}, AMS, Math. Surveys and Monographs, Vol.49, 1997.

\bibitem{So}R.B. Sower: Large deviations for a reaction-diffusion equation with non-Gaussian
perturbations. {\it The Annals of Probability}, 20(1): 504-537, 1992.

\bibitem{St} D.W. Stroock: An Introduction to the Theory of Large Deviations,
Springer-Verlag, New York, 1984.

\bibitem{Zh1}X. Zhang: On Stochastic Evolution Equations with non-Lipschitz Coefficients. Preprint.

\bibitem{Zh2}X. Zhang: A variational representation for random functionals on abstract Wiener spaces. Preprint.

\end{thebibliography}
\end{document}